\newcommand{\References}{references}

\documentclass[a4paper,conference,10pt]{IEEEtran}
\IEEEoverridecommandlockouts

\usepackage{cite}
\usepackage{amsmath,amssymb,amsfonts}
\usepackage{graphicx}
\usepackage{textcomp}
\usepackage{xcolor}
\usepackage{algorithm}
\usepackage[italicComments=true, beginComment=//~]{algpseudocodex}
\def\BibTeX{{\rm B\kern-.05em{\sc i\kern-.025em b}\kern-.08em
    T\kern-.1667em\lower.7ex\hbox{E}\kern-.125emX}}

\usepackage{amsthm}
\usepackage{dsfont}
\usepackage{tikz}
\usetikzlibrary{shapes, arrows, math, calc, shapes.multipart}
\usetikzlibrary{positioning}
\tikzstyle{block} = [draw, fill=white, rectangle, minimum height=3em, minimum width=6em]
\tikzstyle{output} = [coordinate]
\tikzstyle{input} = [coordinate]
\usetikzlibrary{decorations.markings,arrows.meta}
\usepackage{pgfplots}
\pgfplotsset{compat=1.18}

\usepackage{textcomp}   
\def\BibTeX{{\rm B\kern-.05em{\sc i\kern-.025em b}\kern-.08em
    T\kern-.1667em\lower.7ex\hbox{E}\kern-.125emX}}

\usepackage[colorlinks=true]{hyperref}
\usepackage{float}
\usepackage{subcaption}
\usepackage{comment}
\usepackage{enumerate}

\usepackage[normalem]{ulem}
\usepackage{mathtools}
\usepackage{cleveref}

\makeatletter
\newcommand\RedeclareMathOperator{%
  \@ifstar{\def\rmo@s{m}\rmo@redeclare}{\def\rmo@s{o}\rmo@redeclare}%
}
\newcommand\rmo@redeclare[2]{%
  \begingroup \escapechar\m@ne\xdef\@gtempa{{\string#1}}\endgroup
  \expandafter\@ifundefined\@gtempa
     {\@latex@error{\noexpand#1undefined}\@ehc}%
     \relax
  \expandafter\rmo@declmathop\rmo@s{#1}{#2}}
\newcommand\rmo@declmathop[3]{%
  \DeclareRobustCommand{#2}{\qopname\newmcodes@#1{#3}}%
}
\@onlypreamble\RedeclareMathOperator
\makeatother

{\left(\begin{smallmatrix}}%
{\end{smallmatrix}\right)}%

{\left[\begin{smallmatrix}}%
{\end{smallmatrix}\right]}%

\newcommand{\N}{\mathds{N}}
\newcommand{\R}{\mathds{R}}
\newcommand{\Rp}{\R_{\geq0}}

\newcommand{\C}{\mathds{C}}

\newcommand{\fa}{\ \forall \, }

\newcommand{\nl}{\left\|}
\newcommand{\nr}{\right\|}
\newcommand{\cbl}{\left\lbrace }
\newcommand{\cbr}{\right\rbrace }

\newcommand{\Norm}[2][ ]{\nl #2 \nr_{#1}}
\newcommand{\SNorm}[1]{\Norm[\infty]{#1}}

\newcommand{\setdef}[2]{\cbl\ #1\ \left|\ \vphantom{#1} #2\ \right.\cbr}

\newcommand{\GL}{\text{GL}}

\newcommand{\cD}{\mathcal{D}}

\newcommand{\cF}{\mathcal{F}}
\newcommand{\cG}{\mathcal{G}}

\newcommand{\cY}{\mathcal{Y}}

\newcommand{\oT}{\mathbf{T}}

\newcommand{\ds}{d\text{s}}

\DeclareMathOperator*{\rf}{ref}

\newcommand{\ve}{\varepsilon}

\newcommand{\con}{\mathcal{C}}

\RedeclareMathOperator*{\Im}{Im}
\RedeclareMathOperator*{\Re}{Re}
\renewcommand{\phi}{\varphi}

\newcommand{\vp}{\varphi}

\newcommand{\dd}[2][ ]{\tfrac{\text{\normalfont d}#1}{\text{\normalfont d}#2}}

\DeclareMathOperator*{\rank}{rank}

\newcommand{\fM}{L_{\rm max}}
\newcommand{\gM}{\gamma_{\rm max}}
\newcommand{\gm}{\gamma_{\rm min}}

\newcommand{\umax}{u_{\rm max}}

\theoremstyle{plain}
\newtheorem{theorem}{Theorem}[section]
\newtheorem{corollary}[theorem]{Corollary}
\newtheorem{proposition}[theorem]{Proposition}
\newtheorem{assumption}{Assumption}
\newtheorem{definition}[theorem]{Definition}
\newtheorem{remark}[theorem]{Remark}
\newtheorem{lemma}[theorem]{Lemma}
\Crefname{definition}{Definition}{Definitions}
\Crefname{assumption}{Assumption}{Assumptions} 
\Crefname{lemma}{Lemma}{Lemmata}
\Crefname{remark}{Remark}{Remarks}
\Crefname{theorem}{Theorem}{Theorems}
\Crefname{corollary}{Corollary}{Corollaries}
\Crefname{proposition}{Proposition}{Propositions}

\title{\LARGE \bf 
Safe data-driven reference tracking with prescribed performance}
\author{Philipp Schmitz$^{*}$, Lukas Lanza$^{*}$, Karl Worthmann$^{*}$
\thanks{\textsuperscript{$\star$}L.~Lanza and K.~Worthmann gratefully acknowledge funding by the Deutsche Forschungsgemeinschaft (DFG, German Research Foundation) -- {Project-IDs 471539468 and 507037103}.
Philipp Schmitz is grateful for the support from 
the Carl Zeiss Foundation (DeepTurb---Deep Learning in and of Turbulence; project No.\ 2018-02-001).}
\thanks{$^{*}$L.~Lanza, P. Schmitz, and K.~Worthmann are with the Optimization-based Control Group, Institute of Mathematics, Technische Universit\"at Ilmenau, Germany
        {\tt\small \{lukas.lanza, philipp.schmitz, karl.worthmann\}@tu-ilmenau.de}}
}

\begin{document}

\maketitle
\thispagestyle{empty}
\pagestyle{empty}

\begin{abstract}
We study output reference tracking for unknown continuous-time systems with arbitrary relative degree. 
The control objective is to keep the tracking error within predefined time-varying bounds while measurement data is only available at discrete sampling times. 
To achieve the control objective, we propose a two-component controller. One part is a recently developed sampled-data zero-order hold controller, which achieves reference tracking within prescribed error bounds. To further improve the control signal, we explore the system dynamics via input-output data, and include as the second component a data-driven MPC scheme based on Willems et al.'s fundamental lemma. This combination yields significantly improved input signals as illustrated by a numerical example.
\end{abstract}

\section{Introduction}

We study output-reference tracking of unknown continuous-time linear systems with arbitrary relative degree, where the output measurement data is only available at discrete sampling times.
The control task is to ensure that the tracking error evolves strictly within predefined time-varying bounds.
The latter aspect is subject of various results in feedback control, such as prescribed performance control~\cite{bechlioulis2014low}, see also \cite{dimanidis2020output} for a comprehensive literature review; or funnel control~\cite{IlchRyan02b}, cf.~\cite{BergIlch21} for an extensive literature overview.
Since sampled-data control is closely related to digital control, there is plenty of literature available.
In the current context we may pick out~\cite{delchamps1990stabilizing,brockett2000quantized}, where stabilization of linear systems via sampled-data state feedback is studied, and~\cite{chang2019robust}, where $H_\infty$ performance guarantee is achieved by quantized feedback.

Recently, a Zero-order Hold (ZoH) output-feedback controller was proposed~\cite{lanza2023sampleddata}, which achieves reference tracking at the output of an unidentified non-linear system with error guarantees, while only receiving measurement data 
at discrete sampling times. 
However, since only structural assumptions are invoked, the controller 
relies on worst-case estimations and exhibit a potentially unnecessarily large (and heavily oscillating) control signal, 
which is not desirable. 
In this article, we aim to improve the controller performance of the ZoH controller~\cite{lanza2023sampleddata} by combining it with a data-driven Model Predictive Control (MPC) approach. 
The latter is based on the results found by Willems and coauthors in~\cite{WRMDM05}, now known as \emph{Willems et al.'s fundamental lemma}.
This result allows the description of an unknown discrete-time linear time-invariant system in a non-parametric fashion. 
The finite-length input-output trajectories of the system lie in the column space of a suitable Hankel matrix constructed directly from measured data. 
Based on this non-parametric description, it is possible to develop an MPC scheme that uses only input-output data instead of an a priori given model, cf.~\cite{yang2015data,coulson2019data,Berberich20}. 
The fundamental lemma is subject to recent substantial research in the field of data-driven control. In \cite{SchmitzFaulwasserWorthmann22,faulwasser2023behavioral,tudo:pan22a}, 
it was extended to stochastic descriptor systems.
Extensions towards continuous-time and non-linear systems can be found, e.g., in~\cite{lopez2022continuous,rapisarda4370211fundamental,berberich2022linear} and~\cite{Alsalti21,Markovsky21} resp. 
For a 
broader overview,
we refer the reader to~\cite{MarkDorf21} and~\cite{faulwasser2023behavioral}. 

The contribution of the present article is an application of data-based MPC to improve the controller performance for an unknown linear system, while guaranteeing pre-specified time-varying error bounds on the tracking error. 
We partition the set of feasible outputs into a \emph{safe} and a \emph{safety-critical} region. In the safe region, we generate persistently-exciting input signals to explore the system behaviour and, then, set up a learning-based predictive controller (exploitation). In the safety-critical region, we resort to the recently proposed controller~\cite{lanza2023sampleddata} to steer the system into the safe region. 
The basic idea is reminiscent of
~\cite{AlbaAlan16} and follow-up work~\cite{WuAlba19}, where safe and unsafe operating regions were used for economic MPC. 
In comparison to~\cite{lanza2023sampleddata}, we drastically reduce peaks in the control signal. 
Further, we show that the prediction horizon within the MPC scheme can be adaptively and \emph{freely} chosen without jeopardizing (recursive) feasibility and the guaranteed trajectory tracking within the pre-specified error bounds in comparison to existing work, see, e.g.,~\cite{PannWort11,GiseRant13,PannWort14,WortMehr17,Kren18}.

The restriction to linear systems in our presentation is merely of exemplary nature. We emphasize, that for the proposed controller we may also allow for nonlinear continuous-time systems, cf.~\cite{lanza2023sampleddata}. In this case the guaranties on the tracking performance are still valid, even though the significance of the data-based linear surrogate model used in the MPC scheme and, therefore, the quality of the control signal are uncertain. Here, for instance, the mentioned extensions of the fundamental lemma for nonlinear systems have the potential of improving the controller.

This paper is structured in the following way. In \Cref{sec:II} we introduce the examined system class and formulate the control objective. In \Cref{sec:III} we present the controller and its components. To this end, we recap Willems et al.'s fundamental lemma and its leverage in data-driven MPC. Subsequently, we recall the theoretical foundation for the ZoH component. Then, 
an algorithmic implementation of the two-component controller is proposed. In \Cref{sec:IVx} the main result is presented and in \Cref{Sec:Simulation} a numerical example is provided to illustrate the controller performance.

\noindent \textbf{Notation}:
Let $\Rp:=[0,\infty)$. For $x\in\R^n$, the Euclidean norm is denoted by $\Norm{x}$ and ${\Norm{x}}_Q:=\sqrt{x^\top Qx }$ denotes the norm induced by a symmetric positive definite matrix $Q\in\mathbb R^{n\times n}$.
For $p\in\N_0\cup\{\infty\}$ the space of $p$-times continuously differentiable
functions on $V\subset \R^m$ with image in $\R^n$ is denoted by $\con^p(V,\R^n)$ and
$\con(V,\R^n):=\con^0(V,\R^n)$.
Let $L^\infty(I,\R^n)$ be the space of measurable essentially bounded functions 
from an interval $I\subset \R$ into $\R^n$, equipped with the usual norm $\SNorm{\cdot}$. 
$L^\infty_{\text{loc}}(I,\R^n)$ is the space of locally-bounded measurable functions.
$W^{k,\infty}(I,\R^n)$ is the $k$th-order Sobolev space with respect to $L^\infty(I,\R^n)$.
For a finite sequence $(f_k)_{k=0}^{N-1}$ in $\mathbb R^n$ of length $N$ we define the vectorization $f_{[0,N-1]} := \begin{bmatrix}
    f_0^\top & \dots & f_{N-1}^\top
\end{bmatrix}{}^\top \in \mathbb R^{nN}$.
$\GL_m(\R)$ denotes the group of invertible $(m \times m)$ matrices.

\section{Problem formulation}
\label{sec:II}
In this section, we introduce the considered system class and formulate the control objective.

\subsection{System class}
We consider linear time-invariant control systems
\begin{equation} \label{eq:sys}
    \begin{aligned}
        \dot x(t) &= Ax(t) + B u(t), \qquad x(0) = x^0, \\
        y(t) &= C x(t),
    \end{aligned}
\end{equation}
where $x(t) \in \R^n$ is the state at time~$t \geq 0$, $A \in \R^{n \times n}$ is the drift dynamics, $B \in \R^{n \times m}$ is the input distribution matrix, and $C \in \R^{m \times n}$ defines the output matrix. 
Note that input~$u(t)$ and output~$y(t)$ are of the same dimension~$m \in \N$.
We make 
the following assumptions.
\begin{assumption} \label{ass:RelDeg}
    System~\eqref{eq:sys} has strict relative degree~$r \in \N$:
    \begin{equation*}
        \forall \, k = 0,\ldots,r-2 : \, CA^kB = 0, \text{ and } CA^{r-1} B \in \GL_m(\R).
    \end{equation*}
\end{assumption}
\Cref{ass:RelDeg} allows us to find a linear coordinate transformation $U \in \GL_n(\R)$ such that $U x = (y,\dot y,\ldots,y^{(r-1)},\eta)$ represents system~\eqref{eq:sys} in the so-called \emph{Byrnes-Isidori form}~\cite{isidori95}
\begin{subequations} \label{eq:sys_BIF}
    \begin{align}
        y^{(r)}(t) &= \sum_{i=0}^{r-1} L_i y^{(i)}(t) + S \eta(t) + \Gamma u(t), \label{eq:sys_BIF_y}\\
        \dot \eta(t) &= K \eta(t) + P y(t), \label{eq:sys_BIF_eta}
    \end{align}
\end{subequations}
where $(y(0),\dot y(0),\ldots y^{(r-1)}(0),\eta(0)) = U x^0$, and for $i=0,\ldots,r-1$ we have $L_i \in \R^{m \times m}$, $S,P^\top \in \R^{m \times (n-rm)}$, $K \in \R^{(n-rm) \times (n-rm)}$, and $\Gamma := C A^{r-1} B \in \GL_m(\R)$ is the so-called \emph{high-gain matrix}.
The latter is specified in the next assumption.
\begin{assumption} \label{ass:Gam_definite}
    The high-gain matrix $\Gamma = C A^{r-1} B $ has strictly positive definite symmetric part, i.e.
    \begin{equation*}
        \forall \, z \in \R^{m}\setminus \{0\} : \ z^\top (\Gamma+\Gamma^\top) z > 0.
    \end{equation*}
\end{assumption}
Likewise we may allow for strictly negative definite~$\Gamma + \Gamma^\top$ by changing the sign in the feedback law~\eqref{eq:controller_recursive}. \Cref{ass:Gam_definite} ensures the existence of $\gm,\gM > 0$ such that
\begin{equation*} 
    \forall \, z \in \R^m\setminus\{0\} \, : \  \gm \le \frac{z^\top (\Gamma+\Gamma^\top) z}{\|z\|^2} \le \gM .
\end{equation*}
The next assumption concerns the internal dynamics of system~\eqref{eq:sys}, which are given by~\eqref{eq:sys_BIF_eta} in the transformed system~\eqref{eq:sys_BIF}.
\begin{assumption} \label{ass:sys_MIN}
    System~\eqref{eq:sys} is of minimum phase, i.e.,
    \begin{equation*}
        \forall \, \sigma \in \C_- \, : \rank \begin{bmatrix} A-\sigma I & B \\ C & 0 \end{bmatrix} = n+m,
    \end{equation*}
    or equivalently, the matrix~$K$ in~\eqref{eq:sys_BIF_eta} is Hurwitz.
\end{assumption}
Now, we introduce the class of considered systems.
\begin{definition} \label{Def:system-class}
     For $m,r \in \N$ a system~\eqref{eq:sys} belongs to the system class $\Sigma^{m,r}$ if \Cref{ass:RelDeg,ass:Gam_definite,ass:sys_MIN} are satisfied; in this case we write ${(A,B,C)\in \Sigma^{m,r}}$.
\end{definition}
    
Recall from standard theory that given a control $u\in L^\infty(\Rp, \R^m)$ the linear time-invariant system~\eqref{eq:sys_BIF} with initial condition imposed has a unique solution $(y,\eta)$ on $\Rp$.

\subsection{Control objective} \label{Sec:ControlObjective}

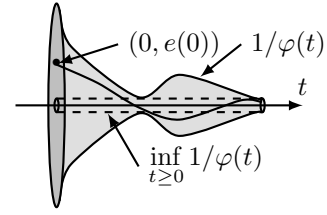
\begin{figure}
  \begin{center}
\begin{tikzpicture}[scale=0.27]
\tikzset{>=latex}
  \filldraw[color=gray!25] plot[smooth] coordinates {(0.15,4.7)(0.7,2.9)(4,0.4)(6,1.5)(9.5,0.4)(10,0.333)(10.01,0.331)(10.041,0.3) (10.041,-0.3)(10.01,-0.331)(10,-0.333)(9.5,-0.4)(6,-1.5)(4,-0.4)(0.7,-2.9)(0.15,-4.7)};
  \draw[thick] plot[smooth] coordinates {(0.15,4.7)(0.7,2.9)(4,0.4)(6,1.5)(9.5,0.4)(10,0.333)(10.01,0.331)(10.041,0.3)};
  \draw[thick] plot[smooth] coordinates {(10.041,-0.3)(10.01,-0.331)(10,-0.333)(9.5,-0.4)(6,-1.5)(4,-0.4)(0.7,-2.9)(0.15,-4.7)};
  \draw[thick,fill=lightgray] (0,0) ellipse (0.4 and 5);
  \draw[thick] (0,0) ellipse (0.1 and 0.333);
  \draw[thick,fill=gray!25] (10.041,0) ellipse (0.1 and 0.333);
  \draw[thick] plot[smooth] coordinates {(0,2)(2,1.1)(4,-0.1)(6,-0.7)(9,0.25)(10,0.15)};
  \draw[thick,->] (-2,0)--(12,0) node[right,above]{\normalsize$t$};
  \draw[thick,dashed](0,0.333)--(10,0.333);
  \draw[thick,dashed](0,-0.333)--(10,-0.333);
  \node [black] at (0,2) {\textbullet};
  \draw[->,thick](4,-3)node[right]{\normalsize$\inf\limits_{t \ge 0} 1/\vp(t)$}--(2.5,-0.4);
  \draw[->,thick](3,3)node[right]{\normalsize$(0,e(0))$}--(0.07,2.07);
  \draw[->,thick](9,3)node[right]{\normalsize $1/\vp(t)$}--(7,1.4);
\end{tikzpicture}
\end{center}
 \vspace*{-2mm}
 \caption{Evolution of the tracking error within the performance margins $\mathcal F_{\vp}$ with boundary $1/\vp(t)$; this figure is based on the figure in~\cite[Fig.~1]{BergLe18a}.}
 \label{Fig:funnel}
 \end{figure}

The control objective is twofold. On the one hand, we aim for an ZoH control input $u\in L^\infty(\Rp,\R^m)$,
\begin{equation*} 
    u(t) \equiv \bar u_k\in\R^m \quad \forall\, t \in [k \cdot \tau, (k+1) \cdot \tau),
    \ k \in \N_0,
\end{equation*}
which achieves that the output~$y$ of system~\eqref{eq:sys} follows a given reference
$y_{\rf}\in W^{r,\infty}(\Rp,\R^{m})$ within pre-specified bounds on the tracking error $e:=y-y_{\rf}$. More precisely, the error $e$ is required to evolve within a given funnel
\begin{align*}
    \cF_\phi= \setdef{(t,e)\in \Rp\times\R^{m}}{\phi(t)\Norm{e} < 1},
\end{align*}
which is determined by a function~$\vp$ being an element of $\cG:=
     \setdef
         {\vp\in W^{1,\infty}(\Rp,\R)}
         {
            \inf_{s \ge 0} \vp(s) > 0
         }$, see \Cref{Fig:funnel}.
Usually, the required tracking performance and, hence, the choice of $\vp$ depend on the particular application. 
This first aspect of the control objective was solved in~\cite{lanza2023sampleddata}.
The second aspect is to significantly improve the controller performance compared to~\cite{lanza2023sampleddata}. This means, we seek a control law, which achieves the aforementioned tracking task, while avoiding the large peaks in the control signal, see the numerical example in.~\cite[Sec.~III]{lanza2023sampleddata}, and~\Cref{fig:ZOH_Controls}.
The latter aspect is achieved by exploring the system dynamics based on input-output data, and then apply data-driven MPC.

\section{Controller structure}
\label{sec:III}
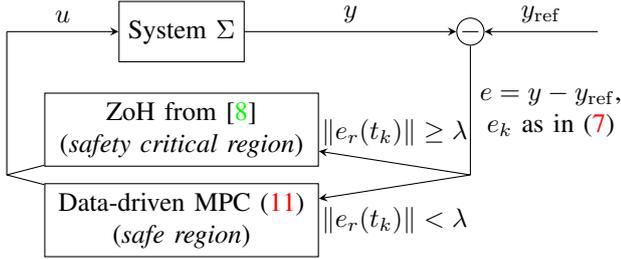
\begin{figure}[H]
  \begin{center}
\begin{tikzpicture} [auto, node distance=2cm,>=stealth, every text node part/.style={align=center}]
\def\hoch{0.8cm};
\def\breit{0cm};
\def\distu{1.3cm};
\def\dista{1.3cm};
\node [block, minimum width = \breit, minimum height = \hoch,] (System) {System~$\Sigma$};
\node [block, minimum width = 3.6cm, minimum height = \hoch, below of=System, node distance = 1.3cm ] (ZoH) {ZoH from~\cite{lanza2023sampleddata} \\ (\emph{safety critical region})};
\node [block, minimum width = 3.6cm, minimum height = \hoch, below of=ZoH, node distance = 1.2cm ] (ZoH+MPC) {Data-driven MPC~\eqref{eq:controller_recursive} \\ (\emph{safe region})};
\node[input, left of=System, node distance = 2.3cm] (u) {};
\coordinate[right of=u, node distance = 0.0cm] (uin) {};
\coordinate[below of=uin, node distance = 1.9cm] (uin_control) ;
\node[circle,,draw=black, fill=white,inner sep=0pt,minimum size=3pt,  right of=System, node distance = 3.8cm] (ref_in) {$-$};
\node[right of=ref_in, node distance = 1.8cm] (ref) {};
\coordinate[below of=ref_in, node distance = 1.9cm] (er);
\draw[->] (u) --  node[above]{$u$} (System);
\draw[->] (System) --  node[above]{$y$} (ref_in);
\draw[->] (ref) --  node[above]{$y_{\rm ref}$} (ref_in);
\draw[->] (er) --  node[yshift=3pt,above]{$\|e_r(t_k)\| \ge \lambda$} (ZoH);
\draw[->] (er) --  node[yshift=-4pt,below]{$\|e_r(t_k)\| < \lambda$} (ZoH+MPC);
\draw[-] (ref_in) -- node[right]{$e = y - y_{\rm ref}$,\\ $e_k$ as in~\eqref{eq:ek}} (er);
\draw[-] (ZoH) -- (uin_control) --  (uin) ;
\draw[-] (ZoH+MPC) -- (uin_control);
\end{tikzpicture}
\end{center}
 \vspace*{-2mm}
 \caption{Schematic structure of the combined controller.}
 \label{Fig:Controller}
 \end{figure}

In this section we introduce the controller. It consists of two components. One component is the sampled-data feedback controller proposed in~\cite{lanza2023sampleddata}, which achieves the control objective defined in~\ref{Sec:ControlObjective}.
We recap this controller in~\Cref{Sec:ZoH}.
The controller~\cite{lanza2023sampleddata} includes an activation threshold~$\lambda$. If the error is below this threshold at sampling~$t_k$, the control signal is determined to be zero for the next interval~$[t_k, t_k+\tau)$. This often results in a large input at the next sampling time, since the system evolves uncontrolled for at least one sampling interval.
However, a closer look to the proof of~\cite[Thm.~2.1]{lanza2023sampleddata} yields that within these intervals of zero input any bounded control can be applied, without jeopardizing feasibility of the controller or violating the error guarantees.
We use this observation for controller design.
Namely, whenever the controller~\cite{lanza2023sampleddata} determines the input to zero, i.e., the system is within a \emph{safe region}, we apply a solution of an optimal control problem (OCP) formulated within a MPC scheme.
In \Cref{Fig:Controller} the controller's structure is depicted schematically.
 For the MPC part of the controller we use a recently developed data-driven framework, cf. \cite{coulson2019data,yang2015data,Berberich20}, introduced in \Cref{Sec:MPC}.
 This MPC algorithm relies on a result from behavioral systems theory in~\cite{WRMDM05}, which is briefly presented in the next \Cref{Sec:FundamentalLemma}.

\subsection{Willems et al.'s fundamental lemma} \label{Sec:FundamentalLemma}

We consider a discrete-time linear surrogate model for the con\-tin\-u\-ous-time system \eqref{eq:sys}, i.e.,
\begin{equation}
\label{eq:discr_sys}
\begin{split}
    x_{k+1} &= \widetilde A x_k + \widetilde B u_k\\
    y_k &= \widetilde C x_k
    \end{split}
\end{equation}
with unknown matrices $\widetilde A \in \mathbb R^{\tilde n\times \tilde n}$, $\widetilde B \in\mathbb R^{\tilde n\times m}$ and $\widetilde C\in\mathbb R^{m\times \tilde n}$. Note that the input and output dimension, respectively, of \eqref{eq:discr_sys} and \eqref{eq:sys} are the same. As the focus below is on input-output trajectories we assume without loss of generality that the system realisation~\eqref{eq:discr_sys} is minimal, i.e.\ controllable and observable. Therefore, it is reasonable to suppose that the state dimension of system~\eqref{eq:discr_sys} is bounded by that of system~\eqref{eq:sys}, that is $\tilde n\leq n$. Instead of $\tilde n$, hereinafter the statements are formulated with respect to the upper bound $n$, which does not change their validity but might come at the expense of higher data demand in applications. 

In the following we recall the notion of persistency of excitation and the fundamental lemma for controllable systems due to Willems et al.~\cite{WRMDM05}. A sequence $(u_k)_{k=0}^{N-1}$ with $u_k\in\mathbb R^m$ is called \emph{persistently exciting of order $L$} for some $L\in\mathbb N$ if the Hankel matrix $\mathcal H_L(u_{[0,N-1]})\in\mathbb R^{mL\times (N-L+1)}$ defined by
\begin{equation} \label{eq:Hankel}
    \mathcal H_L(u_{[0,N-1]}) := \begin{bmatrix}
        u_0 & u_1 & \dots & u_{N-L}\\
        u_1 & u_2 & \dots & u_{N-L+1}\\
        \vdots & \vdots & \ddots &\vdots\\
        u_{L-1} & u_L & \dots & u_{N-1}
    \end{bmatrix},
\end{equation}
has full row rank.
\begin{lemma}[\cite{WRMDM05}]
\label{lem:fl}
    Let $\big((\hat u_k)_{k=0}^{N-1}, (\hat y_k)_{k=0}^{N-1}\big)$ be an input-output trajectory of \eqref{eq:discr_sys} such that the input signal $(\hat u_k)_{k=0}^{N-1}$ is persistently exciting of order $L+n$, where $n$ is the state dimension of system~\eqref{eq:sys}. Then $\big((u_k)_{k=0}^{L-1}, (y_k)_{k=0}^{L-1}\big)$ is an input-output trajectory of \eqref{eq:discr_sys} of length $L$ if and only if there exists $\nu\in \mathbb R^{N-L+1}$ such that
    \begin{equation}
        \begin{bmatrix}
            u_{[0,L-1]}\\
            y_{[0,L-1]}
        \end{bmatrix} =\begin{bmatrix}
            \mathcal H_L(\hat u_{[0,N-1]})\\
            \mathcal H_L(\hat y_{[0,N-1]})
        \end{bmatrix} \nu.
    \end{equation}
\end{lemma}

The fundamental lemma allows a complete description of the finite-length system trajectories in a non-parametric, data-driven manner, without identification of the system matrices, see also the recent survey~\cite{faulwasser2023behavioral} and the references therein for extensions to the descriptor setting including noise.

\subsection{Data-driven MPC component} \label{Sec:MPC}
We aim to track a given sampled output reference signal $(y_{\text{ref},k})_{k=0}^\infty$ while satisfying a control bound  given by some fixed constant $\umax$. Assume that we are given a measured trajectory $\big((\hat u_k)_{k=0}^{N-1}, (\hat y_k)_{k=0}^{N-1}\big)$ of system~\eqref{eq:discr_sys} such that $(\hat u_k)_{k=0}^{N-1}$ is persistently exciting of order $L+2n$. We solve in every discrete time step the optimal control problem
\begin{subequations} \label{eq:OCP_MPC}
\begin{equation}
    \label{eq:opc_dd}
    \operatorname*{minimize}_{\big((u_k)_{k=0}^{L+n-1}, (y_k)_{k=0}^{L+n-1}, \nu\big)} \sum_{k=n}^{L+n-1} \left(\lVert y_k-y_{\text{ref},k}\rVert_Q^2+\lVert u_k\rVert_R^2\right)
\end{equation}
subject to
\begin{align}
\label{eq:OCP_MPC_fl}
    \begin{bmatrix}
            u_{[0,L+n-1]}\\
            y_{[0,L+n-1]}
    \end{bmatrix} &=\begin{bmatrix}
            \mathcal H_{L+n}(\hat u_{[0,N-1]})\\
            \mathcal H_{L+n}(\hat y_{[0,N-1]})
    \end{bmatrix} \nu,\\
    \label{eq:OCP_MPC_init}
    \begin{bmatrix}
            u_{[0,n-1]}\\
            y_{[0,n-1]}
    \end{bmatrix} &= \begin{bmatrix}
            \tilde u_{[0,n-1]}\\
            \tilde y_{[0,n-1]}
    \end{bmatrix},\\
    \Norm{u_k} &\leq u_\text{max},\quad k=n,\dots,L+n-1. \label{eq:ocp_umax}
\end{align}
\end{subequations}
where $\big((\tilde u_k)_{k=0}^{L+n-1}, (\tilde y_k)_{k=0}^{L+n-1}\big)$ is the observed input-output signal from the past to be continued. 
The key difference to standard MPC is that the system dynamics are described by the Hankel matrices in \eqref{eq:OCP_MPC_fl}, cf.\ Lemma~\ref{lem:fl}, instead of a state-space model. The initial condition \eqref{eq:OCP_MPC_init} together with the observability of system~\eqref{eq:discr_sys} guarantees that the latent state is properly aligned, i.e.\ $x_{[0,n-1]} = \tilde x_{[0,n-1]}$ for the state sequences $(x_k)_{k=0}^{L+n-1}$ and $(\tilde x_k)_{k=0}^{n-1}$ corresponding to $\big((u_k)_{k=0}^{L+n-1}, (y_k)_{k=0}^{L+n-1}\big)$ and $\big((\tilde u_k)_{k=0}^{n-1}, (\tilde y_k)_{k=0}^{n-1}\big)$.
The matrices $Q$ and $R$ in \eqref{eq:opc_dd} are assumed to be symmetric positive definite.

\subsection{Zero-order Hold component} \label{Sec:ZoH}
In this section we introduce the ZoH component. To this end, we establish necessary notation and recap some 
results presented in~\cite{lanza2023sampleddata}.
First, we introduce some auxiliary error variables.
Let $\phi\in\cG$, ${y_{\rf}\in W^{r,\infty}(\Rp,\R^m)}$, and
a bijection $\alpha \in\con^1([0,1),[1,\infty))$ be given.
Let $t\geq0$, and set ${{\xi}:=({\xi}_1,\ldots,{\xi}_r)\in\R^{rm}}$.
Then we define
\begin{align} \label{eq:ek}
    e_1(t,{\xi})&:=\phi(t)({\xi}_1-y_{\rf}(t)),\\
    e_{k+1}(t,{\xi})&:=\phi(t)({\xi}_{k+1}-y_{\rf}^{(k)}(t))\!+\!\alpha(\|e_{k}(t,{\xi})\|^2)e_{k}(t,{\xi}), \nonumber
    \end{align}
for $k=1,\ldots,r-1$.
If $\xi(t) = (y(t),\dot y(t),\ldots, y^{(r-1)}(t))$ solves~\eqref{eq:sys_BIF}, then $e_1(t)$ is the tracking error $e(t) = y(t) - y_{\rm ref}(t)$ normalised w.r.t.\ the error boundary~$\vp(t)$.
The bijection $\alpha$ can be suitably chosen, for instance, as $\alpha(s):=1\slash(1-s)$.

The following part contains some technicalities, which are introduced to formulate the controller.
We define the set
\[
        \cD_{t}^{r}:=\setdef
        {\!\! {\xi}\in\R^{rm}}
        { \begin{array}{l}
              \Norm{e_{k}(t,{\xi})}<1,\ k=1,\ldots, r-1, \\
              \|e_r(t,{\xi})\| \le 1
        \end{array} \!\!\! }  .
\]
Using the shorthand notation 
\begin{equation*}
    \chi(y)(t):=(y(t),\dot{y}(t),\ldots,y^{(r-1)}(t))\in\R^{rm}    
\end{equation*}
for $y\in W^{r,\infty}(\Rp,\R^m)$ and $t\in\Rp$,
we define the set of all functions $\zeta\in\con^r([0,\infty),\R^m)$,
which coincide with $\chi(y(0))$, and for which 
$\chi(\zeta)(t)\in\cD^r_t$ on $[0,\delta)$ for $\delta\in(0,\infty]$:
\[
    \cY^r_\delta\!:=\!\setdef
        {\!\!\zeta\in \con^{r-1}([0,\infty),\R^m)\!\!\!}
        {\!\!\!\!\begin{array}{l}
             \chi(\zeta(0))=\chi(y(0)),  \\
             \!\!\fa t\in [0,\delta):\chi(\zeta)(t)\in\cD_t^r\!\!
        \end{array}\!\!\!\!}\!.
\]
This is the set of signals, which evolve ``within the funnel'' to have a simple picture in mind.
For $\zeta \in \cY_\delta^r$ let 
$e_k(t):=e_k(t,\chi(\zeta)(t))$. 
Further, we introduce
the following auxiliary constants $\ve_k,\mu_k$.
Let $\ve_0=0$ and $\bar\gamma_0 :=0$.
We set for $k=1,\ldots,r-1$ successively
\begin{align}
\hat \ve_k \!&\in\! (0,1)  \!\text{ \rm s.t. } 
 \alpha(\hat \ve_k^2) \hat \ve_k \!=\!  \SNorm{\frac{\dot \vp}{\vp} }\!\!\!\!\!\! ( 1\! +\! \alpha(\ve_{k-1}^2) \ve_{k-1})\! +\! 1 \!+\! \bar \gamma_{k-1}, \nonumber \\
    \ve_k \!&:= \max \{ \| e_k(0)\|,  \hat \ve_k\} < 1, \label{eq:ve_mu_gam} \\
    \mu_k \!& := \SNorm{\frac{\dot \vp}{\vp} }\!\!\! ( 1 \! +\!  \alpha(\ve_{k-1}^2) \ve_{k-1} ) \! + \! 1\!+\! \alpha(\ve_k^2) \ve_k  \! + \! \bar \gamma_{k-1} , \nonumber \\
    \bar \gamma_k & := 2 \alpha'(\ve_k^2) \ve_k^2 \mu_k + \alpha(\ve_k^2) \mu_k. \nonumber
\end{align}
For sake of completeness, we recall~\cite[Lem.~1.1]{lanza2023sampleddata} concerning boundedness of the auxiliary error variables~\eqref{eq:ek}.
\begin{lemma}[{\cite[Lem.~1.1]{lanza2023sampleddata}}] \label{Lemma:e_k}
Let $y_{\rf}\in W^{r,\infty}(\Rp,\R^m)$ and $\vp \in \cG$
be given. 
Then, for all $\delta\in(0,\infty]$ and all $\zeta\in \cY^r_\delta$ 
the functions $e_k$ defined in~\eqref{eq:ek} satisfy for the constants $\ve_k>0$, $\mu_k>0$ in~\eqref{eq:ve_mu_gam}
\begin{enumerate}[i)]
    \item $\| e_k(t,\chi(\zeta)(t))\|  \ {\leq}\ \ve_k < 1$,
    \item $\| \dd{t} e_k(t,\chi(\zeta)(t)) \| \ {\leq}\ \mu_k$,
\end{enumerate}
for all $t\in [0,\delta)$ and for all $k=1,\ldots,r-1$.
\end{lemma}

As shown in~\cite{lanza2023sampleddata}, \Cref{Lemma:e_k} allows to infer bounds on the higher derivatives of the signals in~$\cY^r_\delta$, which in particular allows to infer upper bounds on the dynamics of system~\eqref{eq:sys}, if the output~$y(t)$ is ``in the funnel around the reference''.
\begin{corollary}\label{Lemma:DynamicsBounded}
Consider~\eqref{eq:sys} with $(A,B,C)\in\Sigma^{m,r}$.
Let {$y_{\rf} \in W^{r,\infty}(\Rp,\R^m)$}, $\phi\in\cG$, $y^0\in\con^{r-1}(\Rp,\R^m)$ with $\chi(y^0)(0)\in\cD_{0}^r$.
Then, there exists~$\fM \ge 0$, such that for every $\delta\in(0,\infty]$, all $\zeta\in \cY^r_\delta$, and $t\in[0,\delta)$ 
\begin{equation} \label{eq:fmax_gmax_gmin}
\begin{small}
\begin{aligned}
    \fM &\geq \SNorm{ \sum_{i=0}^{r-1} L_i \chi(\zeta)_i \vert_{[0,\delta]} + S \oT(\zeta)\vert_{[0,\delta]} } ,
\end{aligned}
\end{small}
\end{equation}
where $\oT(\zeta)(t) := e^{K t}[0,I_{n-rm}] U x^0 + \int_0^t e^{K(t-s)} P \zeta(s) \ds$ is the solution operator of the internal dynamics~\eqref{eq:sys_BIF_eta}.
\end{corollary}
We omit the proof here, since \Cref{Lemma:DynamicsBounded} is a particular case of~\cite[Lem.~1.2]{lanza2023sampleddata}.
Invoking \Cref{Lemma:e_k} and the constants~\eqref{eq:ve_mu_gam}, and recalling $\eta^0 = [0,I_{n-rm}]U x^0$, we even may calculate explicitly
\begin{equation*}
\begin{small}
\begin{aligned}
  & \SNorm{ \sum_{i=0}^{r-1} L_i \chi(\zeta)_i \vert_{[0,\delta]} + S \oT(\zeta)\vert_{[0,\delta]} } \\
   & \le 
   \sum_{i=0}^{r-1} \| L_i\| \left(\sup_{s \ge 0} \frac{1}{\vp(s)} (1+\ve_i \alpha(\ve_i^2)) + \| y_{\rm ref}^{(i)} \|_\infty \right) \\
   & \quad + \sup_{s \ge 0} \| e^{K s} \| \left( \|P\| ( \sup_{s \ge 0} \frac{1}{\vp(s)} + \| y_{\rm ref} \|_\infty) \int_0^s \| e^{-K \sigma } \| d \text{$\sigma$} \right) \\
   & \quad + \sup_{s \ge 0} \| e^{K s} \| \|\eta^0\| = : \fM.
\end{aligned}
   \end{small}
\end{equation*}
Note that \Cref{Lemma:DynamicsBounded} guarantees an uniform upper bound for the system dynamics, if the output is ``in the funnel around the reference''. 
Although the value~$\fM$ will be used in the control law~\eqref{eq:controller_recursive}, we do not assume knowledge of the system matrices. In particular, the controller is feasible for any $\tilde L_{\rm max} \ge \fM$.
In the context of applications this means that system uncertainties can be estimated roughly, while functioning of the controller is still guaranteed.

\subsection{Controller structure}
We introduce the controller, which achieves the control objective introduced in~\Cref{Sec:ControlObjective}.
With the constants in~\eqref{eq:ve_mu_gam},
we set
\[
        \kappa_0 \!\!:=\!\! \SNorm{\frac{\dot \vp}{\vp}}\!\!\!\!\!\! ( 1 + \alpha(\ve_{r-1}^2) \ve_{r-1} )
                    + \SNorm{\vp}\!\! ( \fM + \|y_{\rf}^{(r)}\|_\infty )+ \bar \gamma_{r-1},
 \]
choose the gain $\beta$, and set the constant~$\kappa_1$ as follows
    \begin{equation*}
        \beta \ge \frac{2 \kappa_0}{ \gm \inf_{s \ge 0} \vp(s) }, \quad 
        \kappa_1 := \kappa_0 +  \SNorm{\vp} \gM \beta.
    \end{equation*}
With the constants defined above, we now provide a uniform bound on the sampling time~$\tau$.
Given an activation threshold ${\lambda \in (0,1)}$, suppose that the sampling time $\tau$ satisfies
    \begin{equation} \label{eq:tau}
0 < \tau \le \min \left\{ \frac{ \kappa_0 }{\kappa_1^2}, \frac{1-\lambda}{\kappa_0 + \|\vp\|_\infty \gM \umax} \right\},
    \end{equation}
where $\umax \ge 0$ is the control bound prescribed in~\eqref{eq:OCP_MPC}.
Let $e_r(t):=e_r(t,(y,\dot{y},\ldots,y^{(r-1)})(t))$, and $t_k := k \tau$.
For~${k \in \N}$ we propose the following controller structure 
    \begin{equation} \label{eq:controller_recursive}
   \forall  t \in  [t_k, t_{k+1})  :  u(t) =
    \begin{cases}
   \quad  u^\star_{\rm MPC}, \!\! &  \|  e_r(t_k)\| < \lambda, \\
         - \beta \tfrac{e_r(t_k)}{\|e_r(t_k)\|^2}  , \!\! &\|  e_r(t_k)\| \ge \lambda
    \end{cases}
    \end{equation}
with activation threshold $ \lambda \in (0,1)$, input gain $\beta > 0$, and $ u^\star_{\rm MPC}$ being the first prospective control action of the solution of the OCP~\eqref{eq:OCP_MPC}. The control scheme based on the control decision~\eqref{eq:controller_recursive} 
is summarized in \Cref{alg:1}.
For $\umax = 0$ the controller~\eqref{eq:controller_recursive} coincides with the controller in~\cite{lanza2023sampleddata}.
Furthermore, the control is uniformly bounded by $\|u(t)\| \le \max \{\beta/\lambda, \umax\}$.
\begin{algorithm}
\caption{Data-driven MPC with error guarantees}\label{alg:cap}
\label{alg:1}
\begin{algorithmic}
\Require $n,L\in\mathbb N$, $\lambda\in (0,1)$
\State $PE\gets \text{false}$;
\While{true}
\State $N\gets N+1$;
\State sample discrete time trajectory $\big((\tilde u_k)_{k=0}^{N-1},(\tilde y_k)_{k=0}^{N}\big)$;
\State calculate $\|e_r\|$;
\If{\textbf{not} $PE$}\Comment{learn the system dynamics}
    \State $\hat u_{N} \gets \tilde u_{N}$, $\hat y_{N} \gets \tilde y_{N}$;
    \If{ $(\hat u_k)_{k=0}^{N-1}$ is p.e.\ of order $L+n$}
        \State $PE\gets \text{true}$;
        \State store $\mathcal H_{L+n}(\hat u_{[0,N-1]})$, $\mathcal H_{L+n}(\hat y_{[0,N-1]})$;
    \EndIf
\EndIf
\If{$\Norm{e_r}<\lambda$}
    \If{PE}
    \Comment{MPC feedback}
        \State $u^\star_\text{MPC}\gets$ solve(OCP~\ref{eq:OCP_MPC});
        \State $u_\text{act} \gets u_\text{MPC}^\star$;
    \Else
    \Comment{random input action}
        \State $u_\text{act} \gets \text{random}$;
    \EndIf
\Else
\Comment{sampled-data feedback}
    \State $u_\text{act} \gets -\beta \frac{e_r}{\Norm{e_r}^2}$;
\EndIf
\State apply $u_\text{act}$ as ZoH input action to the system

\EndWhile
\end{algorithmic}
\end{algorithm}

\begin{remark} \label{Rem:Remark}
    In \Cref{alg:1} the prediction horizon $L$, which is limited by the persistency of excitation order, is fixed. However, the algorithm can be enhanced in the following manner. With increasing time~$t$ more and more system data is available. Hence, a higher persistency of excitation order can be achieved over time, which means that a larger 
    prediction horizon can be used. 
    We illustrate this in a numerical example in \Cref{Sec:Simulation}.
\end{remark}

\section{Main result}
\label{sec:IVx}
In this section we present our main result 
\Cref{Thm:Main}. 
It states that for a given reference~$y_{\rm ref}$ the combined controller~\eqref{eq:controller_recursive} achieves that the output~$y$ of a system~\eqref{eq:sys} tracks the reference with predefined performance.
In particular, \Cref{alg:1} is feasible for all times.
\begin{proposition} \label{Thm:Main}
    Let a reference {$y_{\rf} \in W^{r,\infty}(\Rp,\R^m)$} be given, and choose a funnel function $\vp \in \cG$.
    Consider a system~\eqref{eq:sys} with $(A,B,C) \in \Sigma^{m,r}$.
    Assume
    $\chi(y(0)) \in \cD_{0}^r$, which means that the auxiliary variables~\eqref{eq:ek} (omitting the dependence on $\chi(y) = (y,\ldots,y^{(r-1)})$) satisfy the initial condition
    $\| e_k(0)\| < 1$ for all $k=1,\ldots,r-1$, and $e_r(0) \le 1$;
and, the sampling time~$\tau$ satisfies~\eqref{eq:tau}.
Then the 
controller~\eqref{eq:controller_recursive}
applied to 
system~\eqref{eq:sys} achieves $\| e_k(t)\| < 1$ for all $k=1,\ldots,r-1$, and $\|e_r(t)\| \le 1$ for all $t \ge 0$.
In particular, \Cref{alg:1} is initially and recursively feasible.
Moreover, the tracking error satisfies 
\begin{equation*}
    \forall \, t \ge 0 \, : \ \| y(t) - y_{\rm ref}(t) \| < \frac{1}{\vp(t)}.
\end{equation*}
\end{proposition}
\begin{IEEEproof}
We distinguish the two cases $\umax = 0$ and $\umax > 0$.
If $\umax = 0$, then the controller~\eqref{eq:controller_recursive} coincides with the controller in~\cite{lanza2023sampleddata}.
A closer look at the proof of \cite[Thm.~2.1]{lanza2023sampleddata} shows that if $\|e_r(t_k)\| < \lambda$ any bounded control can be applied for $t \in [t_k, t_k +\tau)$. This requires an adaption of the sampling time~$\tau$, which was made in~\eqref{eq:tau}. 
Since we only incorporate the input constraints $\|u_{\rm MPC}\|_\infty \le \umax$ in~\eqref{eq:ocp_umax} and no output constraints, the OCP~\eqref{eq:OCP_MPC} has a solution for all~$k \in \N$ by construction. 
Therefore, adaption of the sampling time~$\tau$ is sufficient to maintain feasibility of the controller.
\end{IEEEproof}


\section{Numerical example} \label{Sec:Simulation}
\begin{figure}
\begin{center}
\includegraphics[trim=2cm 4cm 5cm 15cm,clip=true,width=4.3cm]{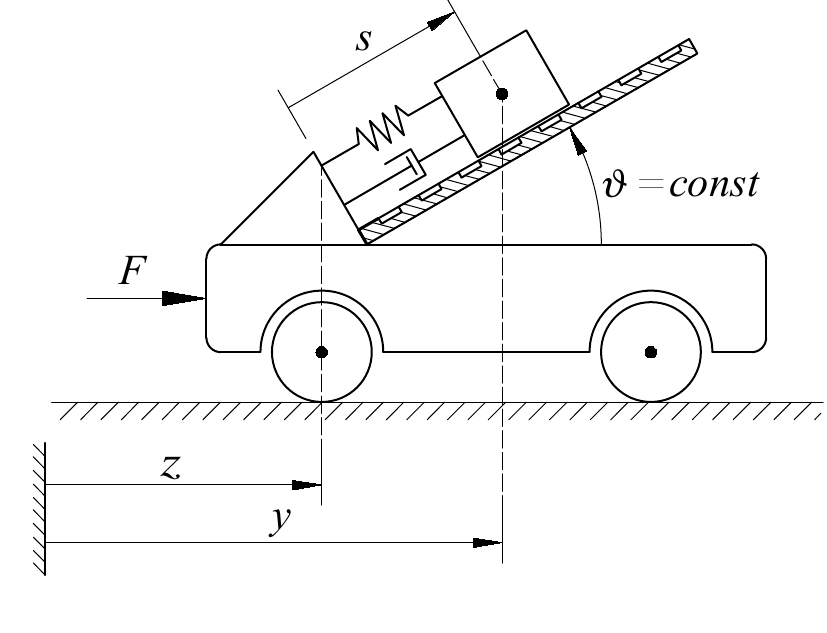}
\end{center}
    \vspace{-5mm}
    \caption{Mass-on-car system. The figure is based on~\cite{SeifBlaj13,BergIlch21}.}
    \label{Fig:Mass-on-a-car}
\end{figure}

We perform a numerical simulation for the mass-on-car system~\cite{SeifBlaj13} to illustrate the proposed control scheme~\eqref{eq:controller_recursive}.
As usual, the mathematical model is given in scaled variables without physical units.
On a car with mass~$m_1$ a ramp is mounted on which a mass~$m_2$ passively slides, see \Cref{Fig:Mass-on-a-car}.
The ramp is inclined by a fixed angle~$\vartheta \in (0,\pi/2)$.
To steer the car, a force $F=u$ can be applied.
A spring-damper combination couples the sliding mass~$m_2$ to the car.
Invoking Newton's mechanical laws, the equations of motion can be derived to be
\begin{subequations}\label{eq:MOC}
\begin{equation}
\begin{small}
    \begin{aligned} 
        \begin{bmatrix}
        m_1 + m_2 & m_2 \cos(\vartheta) \\
        m_2 \cos(\vartheta) & m_2
        \end{bmatrix}
        \! \!
        \begin{pmatrix}
        \ddot z(t) \\ \ddot s(t)
        \end{pmatrix}
        \!+\!
        \begin{pmatrix}
            0 \\
            k s(t) + d s(t)
        \end{pmatrix}
        = 
        \begin{pmatrix}
            u(t) \\ 0
        \end{pmatrix},
    \end{aligned}
\end{small}
\end{equation}
where the horizontal position of the car is~$z$, and the relative position of the sliding mass on the car is~$s$.
The horizontal position of the sliding mass is considered as measured output of the system
\begin{equation}
    y(t) = z(t) + \cos(\vartheta) s(t).
\end{equation}
\end{subequations}
For the numerical simulation we choose the following parameters: inclination angle~$\vartheta = \pi/4$, $m_1 = 1$, $m_2 = 2$, spring constant $k=1$, and damping~$d=1$.
System~\eqref{eq:MOC} can be written as a first order system~\eqref{eq:sys}, and a short calculation yields that for the chosen parameters it satisfies \Cref{ass:RelDeg,ass:sys_MIN,ass:Gam_definite}, with relative degree~$r=2$.
As reference signal we choose $y_{\rf} = 0.4 \sin(\tfrac{\pi}{2}t)$ for $t \in [0,2]$, which means that the mass~$m_2$ on the car is transported from position~$0$ to $0.4$ and back to~$0$, within chosen error boundaries~$\pm 0.15$. 
We start on the reference, i.e., $y(0) = y_{\rm ref}(0)$, $\dot y(0) = \dot y_{\rm ref}(0)$, 
and choose the activation threshold~$\lambda = 0.75$.
As input constraint in the OCP~\eqref{eq:OCP_MPC} we set~$\umax = 20$, and choose $Q=10^2\cdot I$, $R=10^{-4}\cdot I$. Further, we added the regularization term $10^{-6}\cdot\Norm{\nu}^2$ in the cost functional~\eqref{eq:opc_dd}. With these parameters a brief calculation yields $\fM \le 1.4$, $\gm = \gM = 0.25$, and the gain $ \beta \ge 26.98$ is sufficient to guarantee success of the tracking task.
Choosing the smallest~$\beta$ this already gives $\|u_{\rm ZoH+MPC}\|_\infty \le \max\{\beta/\lambda, \umax \}  \le 33.97$.
Moreover, we obtain the sampling time $ \tau \le 4.5 \cdot 10^{-3}$.
We compare the controller~\eqref{eq:controller_recursive} with the controller~\cite{lanza2023sampleddata}, which is~\eqref{eq:controller_recursive} with $u_{\rm MPC}^\star = 0$ (and hence $\umax = 0$ in~\eqref{eq:tau}), as pointed out earlier.
Signals corresponding to the controller~\eqref{eq:controller_recursive} are labeled as $y_{\rm ZoH+MPC}, u_{\rm ZoH+MPC}$ indicating the application of MPC.
Signals corresponding to the controller~\cite{lanza2023sampleddata} are labeled with $ y_{\rm ZoH}, u_{\rm ZoH}$.
Now, we consider two scenarios.
\subsection{Fixed prediction horizon}
\begin{figure}[H]
    \centering
    \includegraphics[scale=0.32]{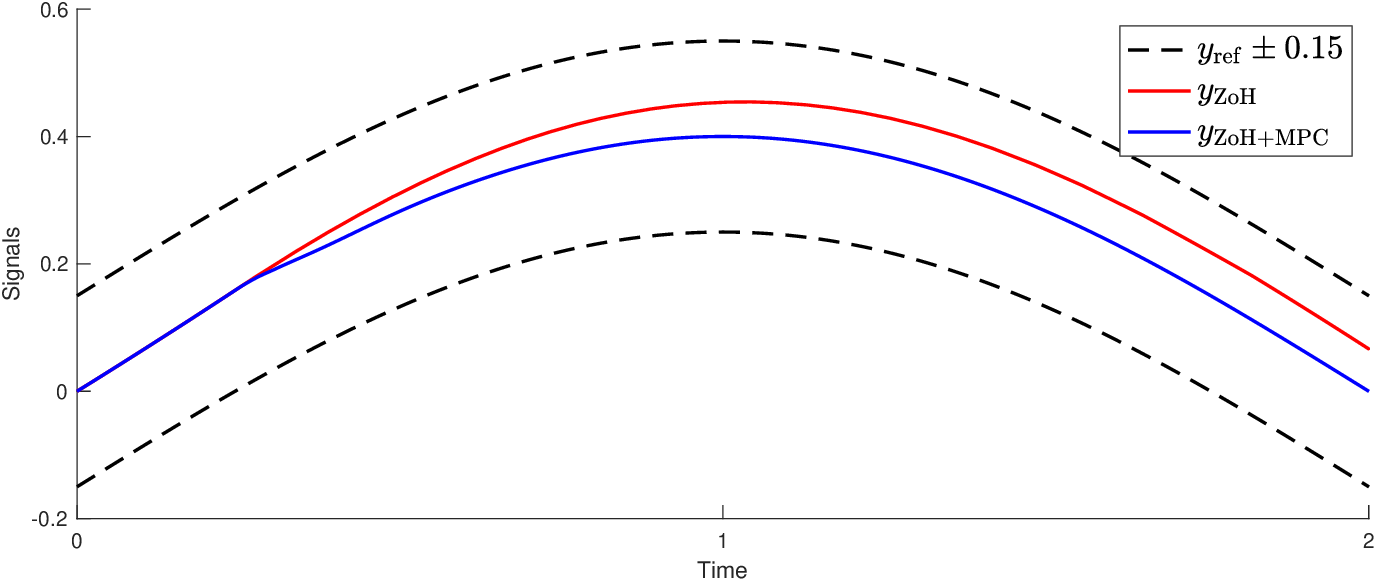}
    \caption{Outputs, reference, and boundaries.}
    \label{fig:ZOH_Errors}
\end{figure}

\begin{figure}[H]
    \centering
    \includegraphics[scale=0.32]{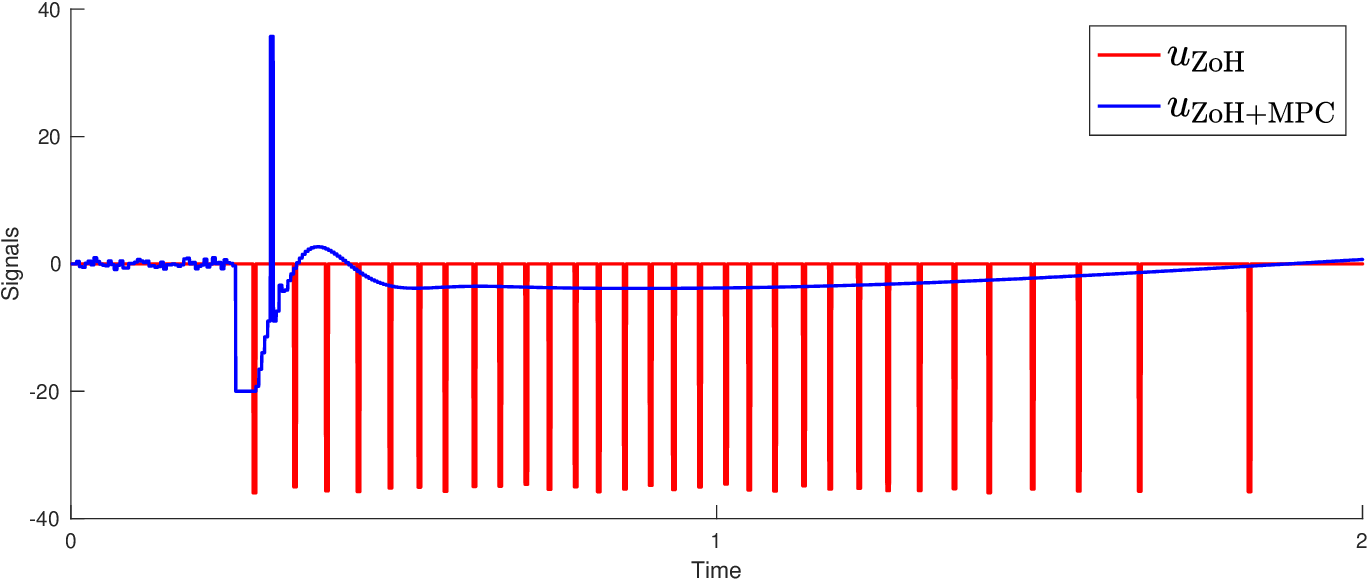}
    \caption{Controls.}
    \label{fig:ZOH_Controls}
\end{figure}

\begin{figure}[H]
    \centering
    \includegraphics[scale=0.32]{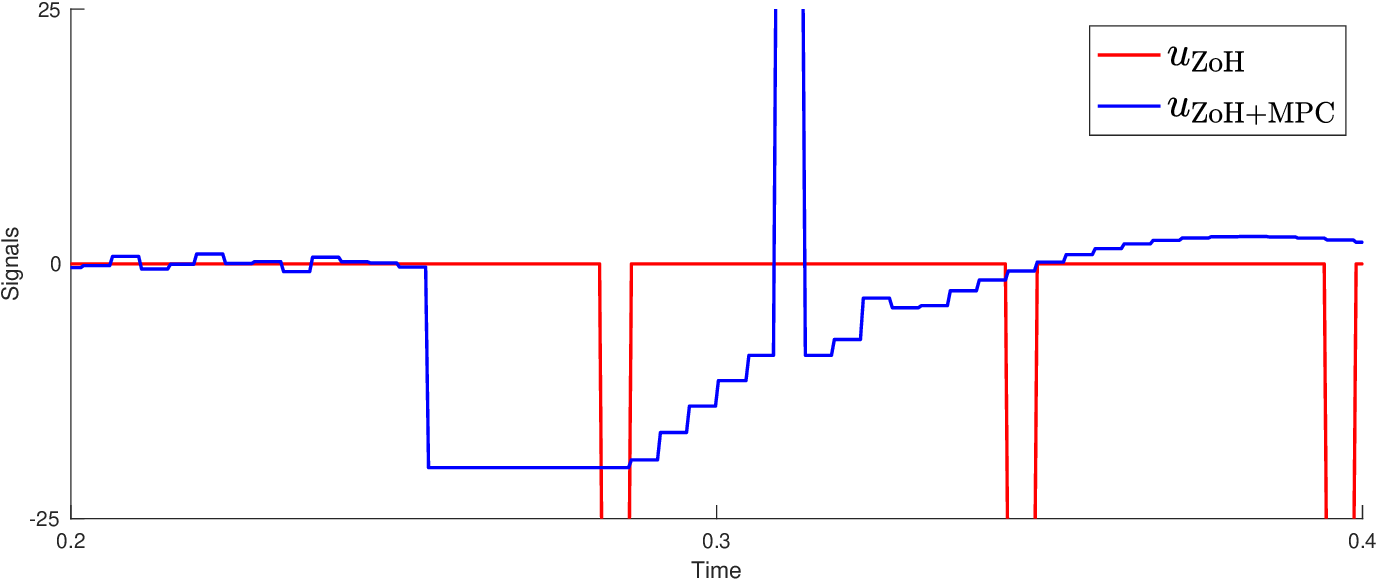}
    \caption{Zoom in controls.}
    \label{fig:ZOH_Controls_zoom}
\end{figure}

\begin{figure}[H]
    \centering
    \includegraphics[scale=0.32]{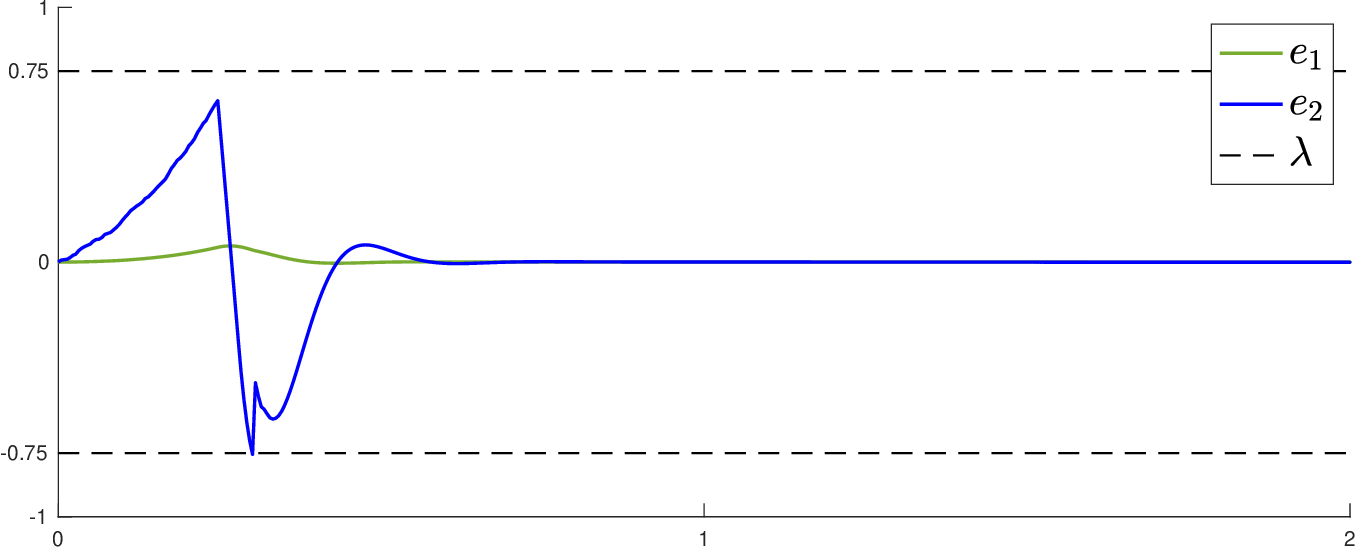}
    \caption{Error variables $e_1, e_2$.}
    \label{fig:e1e2}
\end{figure}

The prediction horizon for MPC is chosen~$L=20$.
In \Cref{fig:ZOH_Errors} the output signals are shown.
Both controllers achieve the tracking task.
However, it can be seen that the controller invoking data-based MPC achieves better tracking performance (the blue line is most of the time very close to the center between the boundaries).
In \Cref{fig:ZOH_Controls} the controls are depicted.
As already observed and discussed in~\cite[Sec.~III]{lanza2023sampleddata}, the control signal~$u_{\rm ZoH}$ consists of separated peaks.
This is due to the incorporation of worst-case estimations in the control law.
Contrary, the control signal~$u_{\rm ZoH+MPC}$ has only one such peak.
This comes from the fact, that the data-based MPC produces control signals~$u^\star_{\rm MPC}$, which are sufficient to keep the error variable $e_2$ below the activation threshold~$\lambda$, see \Cref{fig:e1e2}. In particular, the signal $u_{\rm ZoH+MPC}$ is sufficient to achieve the control objective.
\Cref{fig:ZOH_Controls_zoom} shows a zoom of the control signals.
In the beginning, there is random control in order to generate persistently exciting input signal.
Then, MPC produces a signal, which is saturated by~$\umax$; however, it is still sufficient to keep $e_2$ below~$\lambda$.
Shortly after $t=0.3$ the MPC signal is not sufficient and hence the ZoH signal becomes active, resulting in a large control input, which is applied for one sampling interval. Afterwards, MPC again is sufficient to keep $e_2$ below~$\lambda$.
The evolution of the auxiliary error variables~$e_1, e_2$ is depicted in \Cref{fig:e1e2}.
It can be seen, that only at one sampling instance the case $\| e_2 \| \ge \lambda$ occurs, which then activates the ZoH controller.
%
\subsection{Adaptive prediction horizon}
\begin{figure}[H]
    \centering
    \includegraphics[scale=0.32]{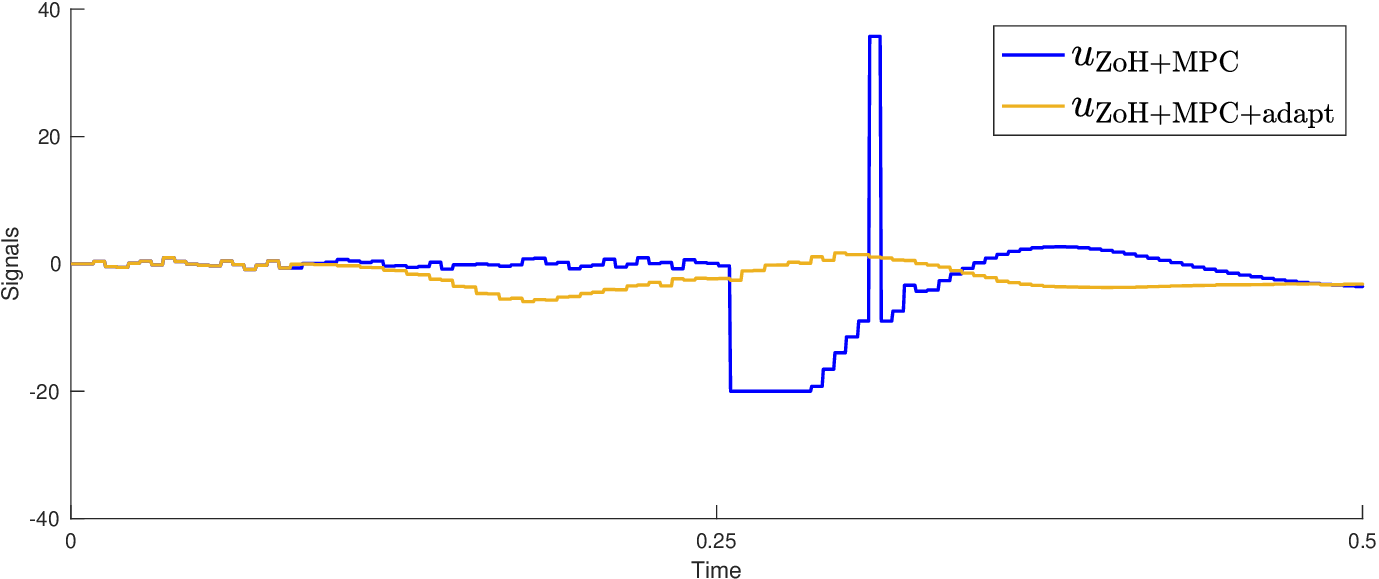}
    \caption{Zoomed in: Controls without and with adapted prediction horizon.}
    \label{fig:in_adapt}
\end{figure}

We adapt the prediction horizon~$L$ in the data-driven MPC, as discussed in \Cref{Rem:Remark}.
Since with increasing time~$t$ more and more system data is available, we update the Hankel matrix~\eqref{eq:Hankel} in every iteration.
Whenever the persistency of excitation order 
increases, we adapt the prediction horizon~$L$.
We start with~$L = 1$ and define an upper limit $L = 50$.
We compare both controller performances, with fixed prediction horizon~$L=20$ (as in the first case), and adaptive horizon.
The tracking performance in both settings is comparable, and looks similar to \Cref{fig:ZOH_Errors}.
We observed, that adapting the prediction horizon in particular affects the beginning of the simulation, 
hence, \Cref{fig:in_adapt} is zoomed in to the beginning.
Since we start with~$L=1$ in the adaptive setting, the MPC is enabled after few iterations, while for fixed horizon L=20 MPC has to wait until the respective persistency of excitation order is reached.
As can be seen from \Cref{fig:in_adapt} this results in a better control signal for adapted horizon.
The control signal with adaptive horizon admits two particularities. 
First, the overall control values are smaller compared to fixed horizon.
Second, with adaptive horizon the MPC is able to completely avoid activation of the ZoH component. Hence, no large peaks occur.

\section{Conclusion and outlook}
We proposed a combined controller, consisting of a sampled-data ZoH feedback component, and a data-driven MPC scheme. 
While the former ensures safety, i.e., guaranteeing the imposed time-varying output constraints, the latter heavily improves the control performance as illustrated in the presented numerical example. In particular, we showed feasibility of the controller in \Cref{Thm:Main}. 
The presented result defines a starting point towards more sophisticated controller designs, which achieve the control objective ``tracking with error guarantees''.
Natural extensions of the presented approach will be the consideration of nonlinear systems, uncertain systems, and incorporation of disturbances.
Moreover, using the recently introduced notion of collective persistency of excitation~\cite{WaarPers20} would further alleviate the applicability of data-based MPC controllers. 

\bibliographystyle{IEEEtran}
\bibliography{\References}

\end{document}